\input ./style/arxiv-general.cfg
\documentclass[MSNbibl,number,citesort,seceqn,dvips]{arxbj}
\makeatletter
   \@ifpackageloaded{graphicx}{}{\usepackage{graphicx}}
\makeatother

\volume{22}
\issue{4}
\pubyear{2016}
\firstpage{2080}
\lastpage{2100}
\doi{10.3150/15-BEJ720}
\docsubty{FLA}

\makeatletter
\renewcommand{\mid}{|}
\newtheorem{theorem0}{Theorem}
\newtheorem{lemma}{Lemma}
\newtheorem{corollary}{Corollary}
\def\bfH{\mathbf{H}}
\def\bfI{\mathbf{I}}

\def\bfX{\mathbf{X}}
\def\bfY{\mathbf{Y}}
\def\bfbeta{\bolds{\beta}}
\def\bftheta{\bolds{\theta}}
\def\bfvareps{\bolds{\varepsilon}}
\newcommand{\go}{\rightarrow}
\newcommand{\vareps}{\varepsilon}
\newcommand{\plim}{\operatorname{plim}}
\makeatother

\begin{document}
\begin{frontmatter}

\title{Consistency of Bayes factor for nonnested model selection when the model dimension~grows}
\runtitle{Bayes factor consistency for nonnested models}

\begin{aug}
\author[A]{\inits{M.}\fnms{Min}~\snm{Wang}\corref{}\thanksref{A}\ead[label=e1]{minwang@mtu.edu}}
 \and
\author[B]{\inits{Y.}\fnms{Yuzo}~\snm{Maruyama}\thanksref{B}\ead[label=e2]{maruyama@csis.u-tokyo.ac.jp}}
\address[A]{Department of Mathematical Sciences, Michigan Technological University, Houghton, MI 49931, USA.\\ \printead{e1}}
\address[B]{Center for Spatial Information Science, University of Tokyo, Bunkyo-ku, Tokyo, 113-0033, Japan.\\ \printead{e2}}
\end{aug}

%
\received{\smonth{12} \syear{2014}}

\begin{abstract}
Zellner's $g$-prior is a popular prior choice for the model selection problems in the context of normal
regression models.
Wang and Sun  [\textit{J. Statist. Plann. Inference} \textbf{147} (2014) 95--105]
recently adopt this prior and put a special hyper-prior for $g$, which results in a closed-form expression of Bayes factor for \textit{nested} linear model comparisons. They have shown that under very general conditions, the Bayes factor is consistent when two competing models are of order $O(n^\tau)$ for $\tau <1$ and for $\tau=1$ is almost consistent except a small inconsistency region around the null hypothesis. In this paper, we study Bayes factor consistency for \textit{nonnested} linear models with a growing number of parameters. Some of the proposed results generalize the ones of the Bayes factor for the case of nested linear models. Specifically, we compare the asymptotic behaviors between the proposed Bayes factor and the intrinsic Bayes factor in the literature.
\end{abstract}

\end{frontmatter}

\section{Introduction}\label{section001}

We reconsider the classical linear regression model
\begin{equation}
\label{model:01} \bfY = {{\mathbf1}_n}\alpha + \bfX_p
\bfbeta_p + \bfvareps,
\end{equation}
where $\bfY = (y_1, \ldots, y_n)'$ is an $n$-vector of responses, $\bfX_p$ is an $n \times p$ design matrix of full column rank, containing all potential predictors, ${\mathbf1}_n$ is an $n \times 1$ vector of ones, $\alpha$ is an unknown intercept, and $\bfbeta_p$ is a $p$-vector of unknown regression coefficients. Throughout the paper, it is assumed that the random error for all models follows the multivariate normal distribution, denoted by $\bfvareps \sim N({\mathbf0}_n, \sigma^2{\bfI}_n)$, where ${\mathbf0}_n$ is an $n \times 1$ vector of zeros, $\sigma^2$ is an unknown positive scalar, and ${\bfI}_n$ is an $n$-dimensional identity matrix. Without loss of generality, we also assume that the columns of $\bfX_p$ have been centered, so that each column has mean zero. 

In the class of linear regression models, we often assume that there is an unknown subset of the important predictors which contributes to the prediction of $\bfY$ or has an impact on the response variable $\bfY$. This is by natural a model selection problem where we would like to select a linear model by identifying the important predictors in this subset. Suppose that we have two such linear regression models $M_j$ and $M_i$, with dimensions $j$ and $i$,
\begin{eqnarray}
\label{model:com01} M_j: \bfY &=& {{\mathbf1}_n}\alpha +
\bfX_j \bfbeta_j + \bfvareps,
\\
\label{model:com02} M_i: \bfY &=& {{\mathbf1}_n}\alpha +
\bfX_i \bfbeta_i + \bfvareps,
\end{eqnarray}
where $\bfX_i$ is an $n \times i$ submatrix of $\bfX_p$ and $\bfbeta_i$ is an $i \times 1$ vector of unknown regression coefficients. As commented by
Kass and Raftery \cite{Kass95}, a natural way to compare the two competing models is the Bayes factor, which has nice model selection consistency properties. Here, consistency means that the true model will be eventually selected if enough data is provided, assuming that the true model exists. Our particular interest in this paper is to study the model selection consistency of Bayes factor when the model dimension grows with the sample size. To be more specific, we consider the following three asymptotic scenarios:
\begin{longlist}[]
\item[Scenario 1.]   $i = O(n^{a_1})$ and $j = O(n^{a_2})$ with $0 \leq a_1 \leq a_2 < 1$.

\item[Scenario 2.]   $i = O(n^{a_1})$ and $j = O(n^{a_2})$ with $0 \leq a_1 < a_2 = 1$.

\item[Scenario 3.]   $i = O(n^{a_1})$ and $j = O(n^{a_2})$ with $a_1 = a_2 = 1$.
\end{longlist}
When the two models $M_i$ and $M_j$ are nested,
Moreno, Gir{\'o}n and Casella
\cite{MoreEliaGirojavi2010} study the consistency of the intrinsic Bayes factor under the three asymptotic scenarios. Later on,
Wang and Sun \cite{WangSun2014} derive an explicit closed-form Bayes factor associated with Zellner's $g$-prior for comparing the two models. They show that under very general conditions, the Bayes factor is consistent when the two models are of order $O(n^\tau)$ for $\tau <1$ and for $\tau=1$ is almost consistent except a small inconsistency region around the null hypothesis. Such a small set of models around the null hypothesis can be characterized in terms of a pseudo-distance between models defined by
Moreno and Gir{\'o}n \cite{MoreGiro2008}. Finally, Wang and Sun \cite{WangSun2014} compare the proposed results with the ones for the intrinsic Bayes factor due to \cite{MoreEliaGirojavi2010}.

It should be noted that $M_i$ and $M_j$ are not necessarily nested in many practical situations. As commented by
Pesaran and Weeks
\cite{PesaWeek1999}, ``\textit{in econometric analysis}, \textit{nonnested models arise naturally when rival economic theories are used to explain the same phenomenon}, \textit{such as unemployment}, \textit{inflation or output growth}.'' In fact, the problem of comparing nonnested models has been studied in a fairly large body of ecomometric and statistical literature from both practical and theoretical viewpoints, dating back to \cite{Hoel1947}. For instance, Cox \cite{Cox1962} develops a likelihood ratio testing procedure and shows that under appropriate conditions, the proposed approach and its variants have well-behaved asymptotic properties. Watnik and Johnson
 \cite{WatnJohn2002} consider the asymptotic behavior of three different testing procedures (the J-test, the JA-test, and the modified Cox test) for the analysis of nonnested linear models under the alternative hypothesis. The interested reader is referred to \cite{WatnJohnBedr2001} and \cite{WatnJohn2002} for detailed descriptions of the three testing procedures. 

Gir{\'o}n et al. \cite{GiroMartMero2006} consider the intrinsic Bayes factor for comparing pairs of nonnested models based on the two different encompassing criteria: \textit{encompassing from above} and \textit{encompassing from below}. Later on, Moreno and Gir{\'o}n \cite{MoreGiro2008} present a comparative analysis of the intrinsic Bayes factor under the two criteria in linear regression models. Recently, Gir{\'o}n et al. \cite{GiroMoreCase2010} study the consistency of the intrinsic Bayes factor for the case of nonnested linear models under the first two asymptotic scenarios above. The latter two papers mainly focus on the consistency of the intrinsic Bayes factor when the model dimension grows with the sample size, whereas under the same asymptotic scenario, the researchers should also be interested in the consistency of Bayes factor based on Zellner's $g$-prior, which is a popular prior choice for the model selection problems in linear regression models. To the best of our knowledge, the latter has just received little attention over the years, even though it is of the utmost importance to address the consistency issue for nonnested models.

In this paper, we investigate Bayes factor consistency  associated with  Zellner's $g$-prior for the problem of comparing nonnested models under the three asymptotic scenarios above. Specifically, we compare the asymptotic results between the proposed Bayes factor and the intrinsic Bayes factor due to \cite{GiroMoreCase2010}. The results show that the asymptotic behaviors of the two Bayes factors are quite comparable in the first two scenarios. It is remarkable that we also study the consistency of the proposed Bayes factor under Scenario $3$, whereas such a scenario is still an open problem for the intrinsic Bayes factor highlighted by Gir{\'o}n et al. \cite{GiroMoreCase2010}.

The remainder of this paper is organized as follows. In Section~\ref{section:02}, we present an explicit closed-form expression of Bayes factor based on the null-based approach. In Section~\ref{section:03}, we address the consistency of Bayes factor for nonnested models under the three asymptotic scenarios. Additionally, we compare the proposed results with the ones of the intrinsic Bayes factor. An application of the results in Section~\ref{section:03} to the ANOVA models is provided in Section~\ref{section:031}. Some concluding remarks are presented in Section~\ref{section:04}, with additional proofs given in the \hyperref[appe]{Appendix}.

\section{Bayes factor} \label{section:02}

Within a Bayesian framework, one of the common ways for the model selection problems is to compare models in terms of their posterior probabilities given by
\begin{equation}
P(M_j\mid \bfY) = \frac{p(M_j)p(\bfY \mid M_j)}{\sum_{i}p(M_i)p(\bfY \mid M_i)} = \frac{p(M_j)\operatorname{BF}[M_{j}: M_b]}{\sum_{i}p(M_i)\operatorname{BF}[M_{i}: M_b]},
\end{equation}
where $p(M_j)$ is the prior probability for model $M_j$ and $p(M_j\mid \bfY)$ is the marginal likelihood of $\bfY$ given $M_j$, and $\operatorname{BF}[M_j: M_b]$ is the \textit{Bayes factor}, which compares each model $M_j$ to the base model $M_b$ and is defined as
\begin{equation}
\label{BF:factor01} \operatorname{BF}[M_j: M_b] =
\frac{p(\bfY \mid M_j)}{p(\bfY \mid M_{b})}.
\end{equation}

The Bayes factor in (\ref{BF:factor01}) depends on the base model $M_b$, which is often chosen arbitrarily in practical situations. There are two common choices for $M_b$: one is \textit{the null-based approach} by using the null model ($M_0$), the other is \textit{the full-based approach} by choosing the full model ($M_F$). This paper focuses on the null-based approach because (i) the null model is commonly used as the base model when using Zellner's $g$-priors in most of the literature \cite{liang2008} and (ii) unlike the full model, the dimension of the null model is independent of the sample size. This is crucial in addressing the consistency of Bayes factor with an increasing model dimension. Accordingly, we compare the reducing model $M_j$ with $M_0$:
\begin{eqnarray}
\label{model:com001} M_j: \bfY &=& {{\mathbf1}_n}\alpha +
\bfX_j \bfbeta_j + \bfvareps,
\\
\label{model:com002} M_0: \bfY &=& {{\mathbf1}_n}\alpha +
\bfvareps.
\end{eqnarray}

Zellner's $g$-prior \cite{Zellner1986} is often to choose the same noninformative priors for the common parameters that appear in both models and to assign Zellner's $g$-prior for others that are only in the larger model. The reasonability of this choice is that if the common parameters are orthogonal (i.e., the expected Fisher information matrix is diagonal) to the new parameters in the larger model, the Bayes factor is quite robust to the choice of the same (even improper) priors for the common parameters; see \cite{KassVaid1992}. Since $\alpha$ and $\sigma^2$ are the common orthogonal parameters in (\ref{model:com001}) and (\ref{model:com002}), we consider the following prior distributions for $(\alpha, \sigma^2, \bfbeta_j)$
\begin{eqnarray}
\label{nice}  M_0: p\bigl(\alpha, \sigma^2\bigr)
&\propto& \frac{1}{\sigma^2},
\nonumber\\[-8pt]\\[-8pt]\nonumber
 M_j: p\bigl(\alpha, \sigma^2, \bfbeta_j
\bigr) &\propto& \frac{1}{\sigma^2} \quad \mbox{and} \quad \bfbeta_j
\mid \sigma^2 \sim N \bigl({\mathbf0}, g\sigma^2 \bigl(
\bfX_j'\bfX_j \bigr)^{-1} \bigr).
\end{eqnarray}
The amount of information in Zellner's $g$-prior is controlled by a scaling factor $g$, and thus the choice of $g$ is quite critical. A nice review of various choices of $g$-priors was provided by Liang et al. \cite{liang2008} and later discussed further by Ley and Steel
\cite{LeySteel2012}. In most of the developments of the $g$-priors, the expression of Bayes factor may not have an analytically tractable form, so numerical approximations will generally be employed, whereas it may not be an easy task for practitioners to choose an appropriate one. In particular, standard approximation, such as Laplace approximation, becomes quite challenging when the number of parameters grows with the sample size.

It is remarkable that Maruyama and George
\cite{MaruGeor2010} propose an explicit closed-form expression of Bayes factor based on combined use of a generalization of Zellner's $g$-prior and the beta-prime prior for $g$:
\begin{equation}
\label{BF:funch1} \pi(g) = \frac{g^b(1 + g)^{-a - b - 2}}{B(a + 1, b + 1)}I_{(0, \infty)}{(g)},
\end{equation}
where $a > -1$, $b > -1$, and $B(\cdot, \cdot)$ is a beta function. Noting that Zellner's $g$-prior is a special case of the generalization of Zellner's $g$-prior in \cite{MaruGeor2010}, we obtain the following result and the proof directly follows Theorem 3.1 of \cite{MaruGeor2010} and is thus omitted for simplicity.

\begin{theorem0}
Under the prior in (\ref{BF:funch1}) with $b = (n-j-1)/2-a-2$, the Bayes factor for comparing $M_j$ and $M_0$ can be simplified as
\begin{equation}
\label{BFequation} \operatorname{BF}[M_j: M_0] =
\frac{\Gamma (j/2+a+1 ) \Gamma ((n-j-1)/2 )}{\Gamma (a+1 )\Gamma ((n-1)/2 )} \bigl(1 - R_j^2 \bigr)^{-(n-j-1)/2+a+1},
\end{equation}
where $R_j^2$ is the usual coefficient of the determination of model $M_j$.
\end{theorem0}

The Bayes factor in (\ref{BFequation}) is very attractive for practitioners because of its explicit expression without integral representation, which is not available for other choices of the hyperparameter~$b$. One may argue that such an expression comes at a certain cost on interpreting the role of the prior for $g$, since this prior depends on both the sample size and the model size through the hyperparameter $b$. It is noteworthy that this type of the prior has been studied in the literature. For example, Bayarri et al. \cite{BayaBergFort2012} propose a truncated version of the beta-prime prior for $g$, such that $g > (n+1)/(j+3)- 1$. A similar type of the prior has also been considered by
Ley and Steel \cite{LeySteel2012}.

At this point, we provide several arguments justifying the specification of the hyperparameters as follows. (i) The choice of $b = (n-j-1)/2-a-2$ yields an implicit $O(n)$ choice of $g$ \cite{MaruGeor2010}, that is, $g = O(n)$, which will prevent the hyper-$g$ prior from asymptomatically dominating the likelihood function; (ii) as the sample size grows, the right tail of the beta-prime prior behaves like $g^{-(a+2)}$, leading to a very fat tail for small values of $a$, an attractive property suggested by
Gustafson, Hossain and MacNab  \cite{GustHoss2006}; (iii) with a choice of $a = -1/2$ and some transformation $\bftheta = (\bfX'\bfX)^{1/2}\bfbeta$, the prior makes the asymptotic tail behavior of
\begin{eqnarray}
p\bigl(\bftheta \mid \sigma^2\bigr) = \int_0^\infty
p\bigl(\bftheta \mid \sigma^2,g\bigr)\pi(g) \,dg
\end{eqnarray}
become the multivariate Cauchy for sufficient large $\bftheta \in R^{p}$, recommended by Zellner \cite{Zellner1986}; (iv)~the resulting Bayes factor in (\ref{BFequation}) enjoys nice theoretical properties and good performances in practical applications; see, for example, \cite{MaruGeor2010,WangSun2014a,WangSun2014}, among others, and (v) when the model dimension $j$ is bounded, the Bayes factor in (\ref{BFequation}) is asymptotically equivalent to the Schwarz approximation.

\begin{theorem0}\label{schwarz:app}
When the model dimension $j$ is fixed, for large sample sizes $n$, the Bayes factor in~(\ref{BFequation}) is equivalent to the Schwarz approximation given by
\begin{equation}
\label{approximation:bf} \operatorname{BF}[M_j: M_0] \approx \exp
\biggl[-\frac{j}{2} \log n - \frac{n}{2}\log \bigl(1 -
R_j^2 \bigr) \biggr].
\end{equation}
\end{theorem0}

\begin{pf}
See the \hyperref[appe]{Appendix}.
\end{pf}

One of the most attractive properties in the Bayesian approaches is the model selection consistency, which means the true model (assuming it exists) will be selected if enough data is provided. This property has been intensively studied under different asymptotic scenarios as the sample size approaches infinity. For example, when the model dimension is fixed, see \cite{liang2008,CaseGeorGir2009,MaruGeor2010,LeySteel2012}, to name just a few. Of particular note is that the consistency of various Bayes factors in the listed references behaves very similarly, because for sufficiently large values of $n$, the intrinsic Bayes factor and Bayes factors associated with mixtures of $g$-priors (e.g., $g = n$ and Zellner--Siow prior) can all be approximated by the Schwarz approximation in (\ref{approximation:bf}); see Theorem 2 of~\cite{MoreGiroCase2014}. Also, we can show that this approximation is valid for the Bayes factor with the hyper-$g$ prior in \cite{liang2008}.


When the model dimension grows with the sample size,
Moreno, Gir{\'o}n and Casella
\cite{MoreEliaGirojavi2010} study the consistency of the intrinsic Bayes factors for comparing nested models, and a generalization of the consistency to \textit{nonnested} models has been addressed by Gir{\'o}n et al. \cite{GiroMoreCase2010}. More recently, Wang and Sun
\cite{WangSun2014} address the consistency of Bayes factor associated with Zellner's $g$-prior for nested models, whereas its consistency for the case of nonnested models is also of the utmost importance. We shall particularly be interested in comparing the asymptotic behaviors between the proposed Bayes factor and the intrinsic Bayes factor under the same asymptotic scenario. The presented results provide researchers a valuable theoretical base for the comparison among nested and nonnested models, which naturally appears in practical situations.

\section{Bayes factor consistency for nonnested linear models} \label{section:03}

In this section, we consider the model selection consistency of Bayes factor  for comparing nonnested models under the three asymptotic scenarios. The Bayes factor in (\ref{BFequation}) may not be directly applied to the problem of comparing nonnested models, whereas we can calculate the Bayes factor between $M_j$ and $M_0$, $\operatorname{BF}[M_j: M_0]$, and the Bayes factor between $M_i$ and $M_0$, $\operatorname{BF}[M_i: M_0]$. Thereafter, the Bayes factor for comparing $M_j$ and $M_i$ can be formulated as
\begin{eqnarray}
\label{nonnest:01} \operatorname{BF}[M_j: M_i] =
\frac{\operatorname{BF}[M_j: M_0]}{\operatorname{BF}[M_i: M_0] }.
\end{eqnarray}
The Bayes factor for comparing $M_j$ and $M_i$ in (\ref{model:com01}) and (\ref{model:com02}) is thus given by
\begin{eqnarray}
\label{nonnest:tw0} \operatorname{BF}[M_j: M_i] =
\frac{\Gamma (j/2+a+1 ) \Gamma ((n-j-1)/2 )}{\Gamma (i/2+a+1 ) \Gamma ((n-i-1)/2 )}\frac{ (1 - R_j^2 )^{-(n-j-1)/2+a+1}}{ (1 - R_i^2 )^{-(n-i-1)/2+a+1}}.
\end{eqnarray}
Let $M_T$ stand for the true model
\[
M_T: \bfY = {{\mathbf1}_n}\alpha + \bfX_T
\bfbeta_T + \bfvareps.
\]
According to \cite{Fern2001}, the Bayes factor is said to be consistent when
\[
\displaystyle\mathop{\plim}_{n\go\infty}\operatorname{BF}[M_j: M_i] = \infty,
\]
if $M_j$ is the true model $M_T$, whereas
\[
\displaystyle\mathop{\plim}_{n\go\infty}\operatorname{BF}[M_j: M_i] = 0,
\]
if $M_i$ is the true model $M_T$, where `plim' stands for convergence in probability and the probability distribution is the sampling distribution under $M_T$. For notational simplicity, let
\[
\delta_{ji} = \frac{1}{\sigma^2}\bfbeta_{j}'
\frac{\bfX_{j}'(\bfI_n-\bfH_{i})\bfX_{p}}{n}\bfbeta_{j},
\]
where $\bfH_{i}= \bfX_{i}(\bfX_{i}'\bfX_{i})^{-1}\bfX_{i}$ with $\bfX_{i}$ being an $n\times i$ submatrix of $\bfX_p$. According to \cite{GiroMoreCase2010}, the value of $\delta_{ji}$ can be viewed as a pseudo-distance between $M_j$ and $M_i$, in which the two models are not necessarily nested. Such a pseudo-distance has the following properties: (i) it is always equal to 0 from any model $M_j$ to itself, that is, $\delta_{jj}=0$; (ii) if $M_i$ is nested in $M_j$,  it is also equal to 0, that is, $\delta_{ij} =0$, and (iii) for any model $M_k$, we have $\delta_{ki} \geq \delta_{kj}$ if $M_i$ is nested in $M_j$. To study the model selection consistency, it is usually assumed that when the sample size approaches infinity, the limiting value of $\delta_{ji}$, denoted by $\delta_{ji}^\ast$, always exists, where
\begin{equation}
\label{delta:fun} \delta_{ji}^\ast = \lim_{n\go\infty}
\frac{1}{\sigma^2}\bfbeta_{j}'\frac{\bfX_{j}'(\bfI_n-\bfH_{i})\bfX_{j}}{n}
\bfbeta_{j}.
\end{equation}

In what follows, let $\lim_{n \go \infty}[M]Z_n$ represent the limit in probability of the random sequence $\{Z_n: n\geq 1\}$ under the assumption that we are sampling from model $M$. We present one useful lemma which is critical for deriving the main theorems in this paper, and the proof of the lemma is directly from Lemma 1 of \cite{GiroMoreCase2010} and is not shown here for simplicity.

\begin{lemma} \label{lemma:00002}
Suppose that we are interested in comparing two models $M_i$ and $M_p$ with dimensions $i$ and $p$, respectively, where $M_i$ is nested in $M_p$. As $n$ approaches infinity, both $i$ and $p$ grow with $n$ as $i = {O}(n^{a_1})$ and $p = O(n^{a_2})$ for $0\leq a_1 \leq a_2 \leq 1$. When sampling from the true model $M_T$,
\begin{longlist}[(iii)]
\item[(i)] if $0 \leq a_1 \leq a_2 < 1$, it follows that
\[
\lim_{n \go \infty}[M_T] \biggl\{\frac{1-R_p^2}{1-R_i^2} \biggr
\} = \frac{1 + \delta_{tp}^\ast}{1 + \delta_{ti}^\ast}.
\]

\item[(ii)] If $0 \leq a_1 < a_2 = 1$, it follows that
\[
\lim_{n \go \infty}[M_T] \biggl\{\frac{1-R_p^2}{1-R_i^2} \biggr
\} = \frac{1 + \delta_{tp}^\ast - 1/r}{1 + \delta_{ti}^\ast},
\]
where $r = \lim_{n\go\infty} n/p >1$.
\item[(iii)] If $a_1 = a_2 = 1$, it follows that
\[
\lim_{n \go \infty}[M_T] \biggl\{\frac{1-R_p^2}{1-R_i^2} \biggr
\} = \frac{1 + \delta_{tp}^\ast - 1/r}{1 + \delta_{ti}^\ast-1/s},
\]
where $r = \lim_{n\go\infty} n/p >1$  and $s = \lim_{n\go\infty} n/i >1$.
\end{longlist}
\end{lemma}

We are now in a position to characterize the consistency of Bayes factor in (\ref{nonnest:tw0}) for comparing nonnested linear models. We begin with Scenario 1, that is, the dimensions of models $M_i$ and $M_j$ are $i = O(n^{a_1})$ and $j=O(n^{a_2})$ with $0  \leq a_1 \leq a_2 <1$, respectively. The following theorem summarizes Bayes factor consistency when either of the two models is the true model.

\begin{theorem0}\label{theorem:002}
Let $M_0$ be the null model nested in both nonnested models $M_i$ and $M_j$, whose dimensions are $i$ and $j$, respectively. Suppose that $i= O(n^{a_1})$ and $j = O(n^{a_2})$ with $0 \leq a_1 \leq a_2 < 1$ and that $\delta_{ij}^\ast >0$ and $\delta_{ji}^\ast >0$. The Bayes factor in (\ref{nonnest:tw0}) is consistent whichever the true model is.
\end{theorem0}

\begin{pf}
See the \hyperref[appe]{Appendix}.
\end{pf}

Under the same asymptotic scenario, Gir{\'o}n et al. \cite{GiroMoreCase2010} also conclude that the intrinsic Bayes factor is consistent whichever the true model is when  $\delta_{ij}^\ast >0$ and $\delta_{ji}^\ast >0$. Such an agreement of the consistency between the two Bayes factors is due to the fact that the dominated term is exactly the same on their asymptotic approximations under Scenario 1. It is noteworthy that Theorem \ref{theorem:002} is also valid for other chosen base model nested in both models $M_i$ and $M_j$, even though the main result of the theorem is derived based on the null-based approach. Moreover, Theorem \ref{theorem:002} can be directly applied to the case in which the dimensions of the two competing models are fixed,\vadjust{\goodbreak} because it can be viewed as a limiting case with both $ \lim_{n\go\infty} n/j$ and $ \lim_{n\go\infty} n/i$ approaching infinity.

\begin{corollary}\label{corollary:1}
Suppose we are interested in comparing two models $M_i$ and $M_j$ with dimensions $i$ and $j$, respectively, and that both dimensions are fixed. The Bayes factor in (\ref{nonnest:tw0}) is consistent under both models provided that $\delta_{ij}^\ast >0$ and $\delta_{ji}^\ast >0$.
\end{corollary}

We now investigate Bayes factor consistency when the dimension of one of the nonnested models is of order $O(n)$. The main results are provided in the following theorem.

\begin{theorem0}\label{theorem:003}
Let $M_0$ be the null model nested in both nonnested models $M_i$ and $M_j$ whose dimensions are $i$ and $j$, respectively. Suppose that $i= O(n^{a_1})$ and $j = O(n^{a_2})$ with $0 \leq a_1 < a_2 = 1$ and that there exists a positive constant $r$ such that $r = \lim_{ n\go\infty}n/j > 1$.
\begin{longlist}[(a)]
\item[(a)] The Bayes factor in (\ref{nonnest:tw0}) is consistent under $M_i$, provided that $\delta_{ij}^\ast > 0$.

\item[(b)] The Bayes factor in (\ref{nonnest:tw0}) is consistent under $M_j$ provided that
\begin{equation}
\label{theorem2:in} \delta_{ji}^\ast \in \bigl(\kappa\bigl(r,
\delta_{j0}^\ast\bigr), \delta_{j0}^\ast
\bigr],
\end{equation}
and $\delta_{j0}^\ast > \delta(r)$, where $\kappa(r, s) =  [r(1+s) ]^{1/r} -1$ and
\begin{equation}
\label{incons:01} \delta(r)= r^{1/(r - 1)} - 1.
\end{equation}
\end{longlist}
\end{theorem0}

\begin{pf}
See the \hyperref[appe]{Appendix}.
\end{pf}

Some of the interesting findings can be drawn from the theorem as follows. First, the lower bound of $\delta_{j0}^\ast $, denoted by $\delta(r)$, is exactly the same as the one in Theorem 2 of \cite{WangSun2014} for comparing nested linear models. Second, Theorem \ref{theorem:003} can be extended to the case of nested model comparisons (i.e., $M_i$ is nested in $M_j$) by assuming that $M_0 = M_i$. Third, the Bayes factor depends on the choice of the base model through the value of $\delta_{j0}^\ast$, and therefore, to enlarger the consistency region in (\ref{theorem2:in}), we need to make $\delta_{j0}^\ast$ be as large as possible. This justifies that the null model $M_0$ would be the best choice as the base model. Fourth, the lower bound of $\delta_{ji}^{\ast}$, denoted by $\kappa(r, \delta_{j0}^\ast)$, is a bounded decreasing function in $r$ and satisfies that for any $\delta_{j0}^\ast>0$,
\[
\lim_{r\go\infty}\kappa\bigl(r, \delta_{j0}^\ast
\bigr) = 0.
\]
Finally, under the same scenario, Gir{\'o}n et al. \cite{GiroMoreCase2010} consider the consistency of the intrinsic Bayes factor and conclude that the intrinsic Bayes factor is consistent under $M_i$ if $\delta_{ij}^\ast>0$ and is consistent under~$M_j$, provided that $\delta_{j0}^\ast > \xi(r)$ with
\begin{eqnarray}
\xi(r) = \frac{r-1}{(r+1)^{(r-1)/r}-1} - 1,
\end{eqnarray}
and
\begin{eqnarray}
\label{theorem2:in1} \delta_{ji}^\ast \in \bigl(\eta\bigl(r,
\delta_{j0}^\ast\bigr), \delta_{j0}^\ast
\bigr],
\end{eqnarray}
where $\eta(r, s) =\frac{r+s}{(1+r)^{(r-1)/r}} - 1$.

\begin{figure}[t]

\includegraphics{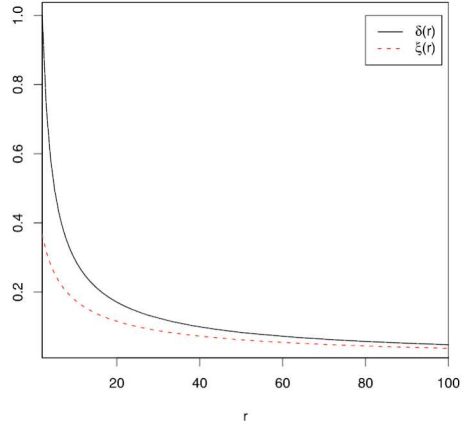}

\caption{The inconsistency region comparisons (below the curves) for the proposed Bayes factor and the intrinsic Bayes factor under Scenario 2.}
\label{fig:0001}
\end{figure}

\begin{figure}[b]
\begin{tabular}{@{}cc@{}}

\includegraphics{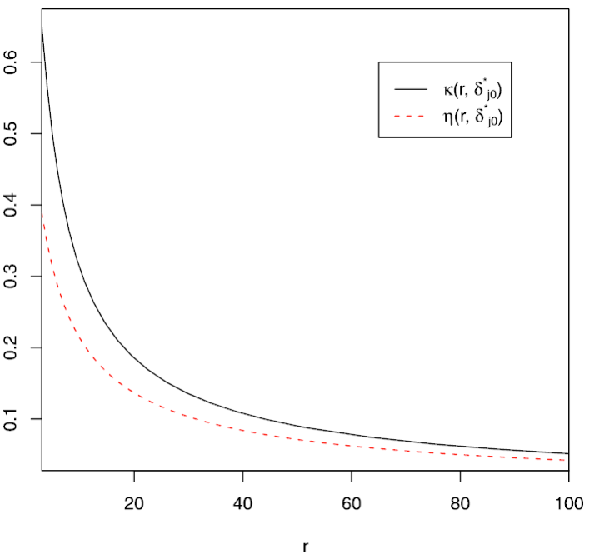}
 & \includegraphics{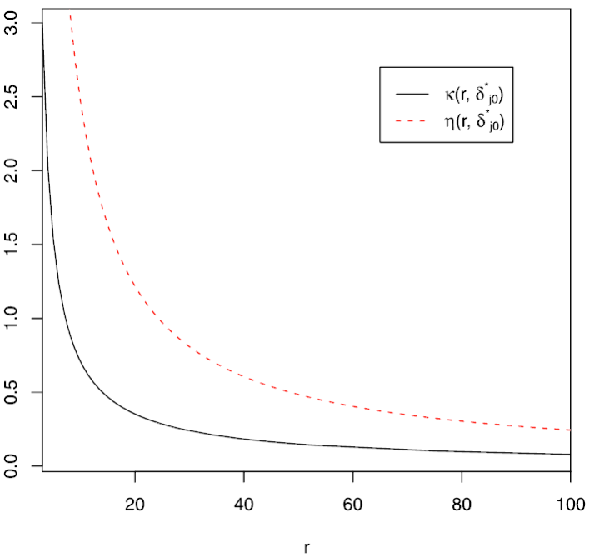}\\
\footnotesize{(a) $\sigma^\ast_{j0}=0.5$} & \footnotesize{(b) $\sigma^\ast_{j0}=20$}
\end{tabular}
\caption{The lower bounds of the consistency regions in (\protect\ref{theorem2:in})  and (\protect\ref{theorem2:in1}) with different limiting values of $\delta_{j0}$ under Scenario~2.}\label{fig:0002}\label{fig:00021}\label{fig:00022}
\end{figure}

It is interesting to observe that the asymptotic behaviors of the two Bayes factors depend on the pseudo-distance between models $\delta_{ji}^\ast$ bounded by $\delta_{j0}^\ast$. Figure~\ref{fig:0001} shows that the upper bounds of their inconsistency regions tend to each other as $r$ increases. Moreover, Figure~\ref{fig:0002} provides their lower bounds with different values of $\delta_{j0}^\ast$. When $\delta_{j0}^\ast$ is small, the consistency region of the proposed Bayes factor is included by the one of the intrinsic Bayes factor, whereas the difference between the two regions is small; see Figure~\ref{fig:00021}(a). However, when $\delta_{j0}^\ast$ gets larger, the consistency region of the proposed Bayes factor will contain the one of the intrinsic Bayes factor, whereas the difference between the two regions becomes significantly as $\delta_{j0}^\ast$ increases; see Figure~\ref{fig:00022}(b). Thus, we may conclude that as $\delta_{j0}^\ast$ increases,\vspace*{1pt} the proposed Bayes factor outperforms the intrinsic Bayes factor from a theoretical viewpoint.

It deserves mentioning that the existence of an inconsistency region around the null hypothesis is quite reasonable from a practical point of view, because the nontrue smaller model $M_i$ is parsimonious under large-$p$ situation and is generally selected when conducting model selection, if the true larger model $M_j$ is not so distinguishable from $M_i$. From the prediction view of point, Maruyama \cite{Maruyama2009} has demonstrated the reasonability of the inconsistency region for the one-way fixed-effect ANOVA model, which could be viewed as a special case of the classical linear models in (\ref{model:01}) after some reparameterization. A theoretical justification of this line of thought for a more general model is still under investigation and will be reported elsewhere.

The first two theorems mainly focus on the consistency of Bayes factor for the case in which at least one model is of order $O(n^\alpha)$ for $\alpha<1$. It is worthy of investigating the consistency issue for the case where both models are of order $O(n)$: the growth rates of the two model dimensions are as fast as $n$. Such a scenario remians an open problem for the intrinsic Bayes factor commented by Gir{\'o}n et al. \cite{GiroMoreCase2010}. We summarize the consistency of the proposed Bayes factor under this scenario in the following theorem.

\begin{theorem0}\label{theorem:004}
Let $M_0$ be the null model nested in both nonnested models $M_i$ and $M_j$ with dimensions $i= O(n)$ and $j = O(n)$, respectively. Suppose that there exist positive constants $r$ and $s$ such that $r = \lim_{ n\go\infty}n/j > 1$ and $s = \lim_{ n\go\infty}n/i > 1$. Without loss of generality, we assume that $r \leq s$.
\begin{longlist}[(a)]
\item[(a)] The Bayes factor in (\ref{nonnest:tw0}) is consistent under $M_i$ provided that
\begin{equation}
\label{inthm3:02} \delta_{ij}^\ast \in \biggl(\frac{r-1}{r}
\biggl\{ \biggl[\frac{s^{1/s}}{r^{1/r}}\bigl(1+\delta_{i0}^\ast
\bigr)^{1/s - 1/r} \biggr]^{r/(r-1)} - 1 \biggr\}, \delta_{i0}^\ast
 \biggr],
\end{equation}
and that $\delta_{i0}^\ast > 0$ satisfying
\begin{equation}
\label{inthm3:01} \biggl(1 + \frac{\delta_{i0}^\ast}{1-1/r} \biggr)^{1-1/r} >
\frac{(1/r)^{1/r}}{(1/s)^{1/s}} \bigl(1+\delta_{i0}^\ast
\bigr)^{1/s-1/r}.
\end{equation}
\item[(b)] The Bayes factor in (\ref{nonnest:tw0}) is consistent under $M_j$ provided that
\begin{equation}
\label{lower:01} \delta_{ji}^\ast \in \bigl(\phi \bigl(r, s,
\delta_{j0}^\ast \bigr), \delta_{j0}^\ast
\bigr],
\end{equation}
where
\[
\phi(a, b, c) =\frac{b-1}{b} \biggl[\frac{a^{1/a}}{b^{1/b}} (1 + c
)^{1/a-1/b} - 1 \biggr]^{b/(b-1)},
\]
and that $\delta_{j0}^\ast > 0$ satisfying
\begin{equation}
\label{inthm3:03} \biggl(1 + \frac{\delta_{j0}^\ast}{1-1/s} \biggr)^{1-1/s}>
\frac{r^{1/r}}{s^{1/s}}\bigl(1+\delta_{j0}^\ast
\bigr)^{1/r-1/s}.
\end{equation}
\end{longlist}
\end{theorem0}

\begin{pf}
See the \hyperref[appe]{Appendix}.
\end{pf}

Unlike the first two asymptotic scenarios, Theorem \ref{theorem:004}(a) shows that under Scenario 3, there exists an inconsistency region around the alternative hypothesis when $M_i$ is true and that the consistency under $M_i$ depends on the chosen base model $M_0$ through the distance $\delta_{i0}^\ast$ only. The existence of the inconsistency region is quite reasonable because there are many candidates to be the base model, which could have a dimension of order $O(n^{a_1})$ with $a_1 \leq 1$. In particular, we observe that the inconsistency region disappears for the case in which $r = s$. This is also very understandable, because with the same growth rates, the parsimonious model is typically preferred in terms of model selection. Furthermore, it can be easily shown that the inequality in (\ref{inthm3:01}) and the lower bound of the consistency region in (\ref{inthm3:02}) are both valid for any $\delta_{i0}^\ast > 0$ if $s^{1/s} \leq r^{1/r}$, indicating that for any $\delta_{i0}^\ast > 0$, the inconsistency region disappears whenever $s \geq r \geq e \approx 2.718$. In order to enlarger the consistency region in (\ref{inthm3:02}), we need to choose a base model to maximize the distance $\delta_{i0}^\ast$. Finally, when $s$ tends to infinity, the inconsistency region disappears for any $\delta_{i0}^\ast > 0$ and $r > 1$, 
which shows that Theorem \ref{theorem:004}(a) just reduces to Theorem~\ref{theorem:003}(a).

\begin{table}[b]
\tabcolsep=0pt
\caption{The consistency regions of the Bayes factor in (\protect\ref{nonnest:tw0}) and the intrinsic Bayes factor due to \cite{GiroMoreCase2010} for different choices of $a_1$ and $a_2$}
\label{table21:ch2}
\begin{tabular*}{\tablewidth}{@{\extracolsep{\fill}}@{}lll@{}}
\hline
Rate of divergence        & The proposed Bayes factor                             & The intrinsic Bayes factor \\
\hline
 $0 < a_1 = a_2 = 1$       & $M_j$: $\delta_{j0}^\ast > \psi(r)$ and $\delta_{ji}^\ast  \in  (\phi (r, s, \delta_{j0}^\ast ), \delta_{j0}^\ast  ]$          & $M_j$: unknown  \\[2pt]
 $0 \leq a_1 < a_2 = 1$    & $M_j$: $\delta_{j0}^\ast > \delta(r)$ and $\delta_{ji}^\ast \in  (\kappa(r, \delta_{j0}^\ast), \delta_{j0}^\ast  ]$          & $M_j$: $\delta_{j0}^\ast > \xi(r)$ and $\delta_{ji}^\ast \in  (\eta(r, \delta_{j0}^\ast), \delta_{j0}^\ast ]$        \\[2pt]
 $0 \leq a_1 \leq a_2 < 1$ & $M_j$: $\delta_{ij}^\ast > 0$ and $\delta_{ji}^\ast > 0$& $M_j$: $\delta_{ij}^\ast > 0$ and $\delta_{ji}^\ast > 0$                 \\
\hline
\end{tabular*}
\end{table}

Theorem \ref{theorem:004}(b) shows that the consistency region under $M_j$ depends on the chosen base model through $\delta_{j0}^\ast$ only. Thus,\vspace*{1pt} the base model should be chosen as small as possible to maximize the value of $\delta_{j0}^\ast$. Note that when $r = s$, the inconsistency region disappears under $M_j$. Also, if the rate of growth of $M_i$ is smaller than that of $M_j$ (i.e., $s$ tends to infinity), then with $\lim_{s\go\infty} s^{1/s} = 1$, the inequality in (\ref{inthm3:03}) turns to be
\begin{equation}
\delta_{j0}^\ast > r^{1/(r - 1)} - 1 =\delta(r),
\end{equation}
which becomes inequality in (\ref{incons:01}) in Theorem \ref{theorem:003}, and the lower bound in (\ref{lower:01}) is
\begin{eqnarray*}
\lim_{s\go\infty}\phi \bigl(r, s, \delta_{j0}^\ast
\bigr) &=& \lim_{s\go\infty}\frac{s-1}{s} \biggl[
\frac{r^{1/r}}{s^{1/s}} \bigl(1 + \delta_{j0}^\ast
\bigr)^{1/r-1/s} - 1 \biggr]^{s/(s-1)}
\\
&=& \bigl[r\bigl(1+\delta_{j0}^\ast\bigr)
\bigr]^{1/r} -1 = \kappa\bigl(r, \delta_{j0}^\ast
\bigr).
\end{eqnarray*}
This illustrates that Theorem \ref{theorem:003}(b) is just a special of Theorem \ref{theorem:004}(b) when $s$ approaches infinity. We may thus conclude that when $s$ tends to infinity, Theorem \ref{theorem:004} reduces to Theorem \ref{theorem:003}.

We have compared the consistency of the proposed Bayes factor with the one of the intrinsic Bayes factor due to \cite{GiroMoreCase2010} under the first two asymptotic scenarios above. A brief summary of comparisons between the two Bayes factors is presented in Table~\ref{table21:ch2}. We observe that the consistency results presented here are similar to the ones for the intrinsic Bayes factor studied by Gir{\'o}n et al. \cite{GiroMoreCase2010}. The similarity occurs, mainly because the asymptotic behaviors of the two Bayes factors depend on a limiting value of $(1-R_j^2)/(1-R_i^2)$ summarized in Lemma \ref{lemma:00002}. The consistency of the intrinsic Bayes factor is still an open problem under Scenario 3. We presume that under Scenario 3, the consistency of the intrinsic Bayes factor also behaves similarly with the one of the proposed Bayes factor, but some further investigation about this presumption is required.

\section{Application} \label{section:031}

It is well known that the ANalysis Of VAriance (ANOVA) models are extremely important in exploratory and confirmatory data analysis in various fields, including agriculture, biology, ecology, and psychology studies. One major difference between the ANOVA models and the classical linear model is that the matrix $[{\mathbf1}_n, {\bfX}_p]$ does not necessarily have full column rank in ANOVA setting. Some constraints are thus required for making the model be identifiable. Here, under the sum-to-zero constraint \cite{Fuji1993}, the ANOVA model with constraints for uniqueness can be reparameterized into the classical linear model without constraints; see \cite{WetzGras2012}.

As an illustration, Maruyama \cite{Maruyama2009} and Wang and Sun
\cite{WangSun2012} reparameterize the ANOVA models with the sum-to-zero constraint into the classical linear model in (\ref{model:01}). Thereafter, based on Zellner's $g$-prior with the beta-prime prior for $g$, they obtain an explicit closed-form Bayes factor, which can be treated as a special case of the Bayes factor in (\ref{BFequation}). Consequently, the asymptotic results of the proposed Bayes factor can be easily applied to various ANOVA models. The application to the one-way ANOVA model is straightforward and is thus omitted here for simplicity. In this section, we mainly consider the results for the two-way balanced ANOVA model with the same number of observations per cell. It deserves mentioning that the results can also be generalized to cover the unbalanced case.

Consider a factorial design with two treatment factors $A$ and $B$ having $p$ and $q$ levels, respectively, with a total of $pq$ factorial cells. Suppose $y_{ijl}$ is the $l$th observation in the $(i, j)$th cell defined by the $i$th level of $A$ and the $j$th level of $B$, satisfying the following model
\begin{eqnarray}
\label{anova:01} y_{ijl} = \mu + \alpha_i +
\beta_j + \gamma_{ij} + \vareps_{ijl}, \qquad
\vareps_{ijl} \sim N\bigl(0, \sigma^2\bigr),
\end{eqnarray}
for $i=1, \ldots, p$, $j = 1, \ldots, q$, and $l =1, \ldots, r$. The number of parameters is $pqr$. We shall be interested in the following five submodels:
\begin{longlist}
\item[$M_0$:]   No effect of $A$ and no effect of $B$, that is, $\alpha_i = 0, \beta_j = 0$, and $\gamma_{ij} = 0$ for all $i$ and $j$.

\item[$M_1$:]   Only effect of $A$, that is, $\beta_j=0$ and $\gamma_{ij} = 0$ for all $i$ and $j$.

\item[$M_2$:]   Only effect of $B$, that is, $\alpha_i=0$ and $\gamma_{ij} = 0$ for all $i$ and $j$.

\item[$M_3$:]   The additive model (without interaction), that is, $\gamma_{ij} = 0$  for all $i$ and $j$.

\item[$M_4$:]   The full model (with interaction).
\end{longlist}

By using the sum-to-zero constraint, Maruyama derives an explicit closed-form Bayes factor associated with Zellner's g-prior for the regression coefficients of the reparameterized model (i.e., equation (4.7) of \cite{Maruyama2009}) and the beta-prime distribution for the scaling factor $g$. Moreover, Maruyama studies the consistency of Bayes factor under different asymptotic scenarios. When both $p$ and $q$ approach infinity and $r$ is fixed, Maruyama concludes that the Bayes factor is consistent except under the full model $M_4$, and that when sampling from $M_4$, the Bayes factor is consistent only if
\begin{eqnarray}
\delta_{43}^\ast > H\bigl(r, \delta_{10}^\ast+
\delta_{20}^\ast\bigr),
\end{eqnarray}
where\vspace*{1pt} $\delta_{ji}^\ast$ is equal to the limit of the sum of squares of the differences between the coefficients of model $M_i$ and the coefficients of model $M_j$ as $n$ tends to infinity, and $H(r,c)$ with positive $c$ is the (unique) positive solution of
\begin{eqnarray}
\label{maru:fun} \frac{(x+1)^r}{r} - (x+1) - c = 0.
\end{eqnarray}
Such an inconsistency region occurs due to the model comparison between $M_4$ and $M_3$. Of particular note is that when comparing $M_4$ and $M_3$, we are in the case of Theorem \ref{theorem:003} with $a_2=1$ and that any null hypothesis will result in a model $M_i$ with a reduced set of parameters that will satisfy $a_1 < a_2$ of Theorem \ref{theorem:003}. Consequently, when sampling from the full model $M_4$, the Bayes factor in (\ref{nonnest:tw0}) is consistent only if $\delta_{4i}^\ast \leq \delta_{40}^\ast$ and
\begin{equation}
\label{anova:mycase} \delta_{4i}^\ast > \bigl[r\bigl(1+
\delta_{40}^\ast\bigr) \bigr]^{1/r} - 1.
\end{equation}
When comparing models $M_{4}$ and $M_3$, the consistency region in (\ref{anova:mycase}) becomes
\[
\delta_{43}^\ast > \bigl[r\bigl(1+\delta_{10}^\ast
+ \delta_{20}^\ast + \delta_{43}^\ast
\bigr) \bigr]^{1/r} - 1,
\]
which is equivalent to
\begin{eqnarray}
\frac{(\delta_{43}^\ast+1)^r}{r} - \bigl(\delta_{43}^\ast+1\bigr) - \bigl(
\delta_{10}^\ast+ \delta_{20}^\ast\bigr) =
0.
\end{eqnarray}
This is exactly coincident with equation (\ref{maru:fun}) provided by Maruyama \cite{Maruyama2009}. It deserves mentioning that an extension of the results of the preceding section to higher-order designs is straightforward.

\section{Concluding remarks} \label{section:04}

In this paper, we have investigated the consistency of Bayes factor for nonnested linear models for the case in which the model dimension grows with the sample size. It has been shown that in some cases, the proposed Bayes factor is consistent whichever the true model is, and that in others, the consistency depends on the pseudo-distance between the larger model and the base model. Specifically, the pseudo-distance can be used to characterize the inconsistency region of Bayes factor. By comparing the consistency issues between the proposed Bayes factor and the intrinsic Bayes factor, we observe that the asymptotic results presented here are similar to the ones for the intrinsic Bayes factor. It would be interesting to see the finite sample performance of the two Bayes factors, which is currently under investigation and will be reported elsewhere.

The consistency of Bayes factor further indicates that besides the three commonly used families of hyper-$g$ priors in \cite{liang2008}, the beta-prime prior is also a good candidate for the scaling factor $g$ in Zellner's $g$-prior. Such a comment has also been claimed by Wang and Sun \cite{WangSun2014} when studying Bayes factor consistency for nested linear models with a growing number of parameters. From a theoretical point of view, we may conclude that like the intrinsic Bayes factor, the proposed Bayes factor should also serve as a powerful tool for model selection in the class of normal regression models due to its comparable asymptotic performance.

It is worth investigating the issues of consistency of Bayes factor based on mixtures of \mbox{$g$-}priors due to \cite{liang2008} under the three asymptotic scenarios. However, in most of the developments of the $g$-priors, the expression of Bayes factor may not have an analytically tractable form, and some efficient approximations are required. Standard approximation technique, such as Laplace approximation, becomes quite challenging when the number of parameters grows with the sample size, because the error in approximations needs to be uniformly small over the class of all possible models. Such a situation has also been encountered by
Berger, Ghosh and Mukhopadhyay
\cite{BergGhosMukh2003} when studying the ANOVA models. We plan to address these issues in our future work.

Finally, it deserves mentioning that we mainly address Bayes factor consistency based on a special choice of the hyperparameter $b$ in the beta-prime prior, which results in an explicit closed-form expression of Bayes factor. In an ongoing project, we investigate the effects of $b$ on the consistency of Bayes factor, especially for the case when $b$ does not actually depend on $n$.

\begin{appendix}\label{appe}
\section*{Appendix}

It is well known that the asymptotic approximation of the gamma function, given by Stirling's formula, can be approximated by
\begin{equation}
\label{stiring:1} \Gamma{(\gamma_{1}x +\gamma_{2})} \approx
\sqrt{2\pi}e^{-\gamma_{1}x} (\gamma_{1}x )^{\gamma_{1}x + \gamma_{2} - 1/2},
\end{equation}
when $x$ is sufficiently large. Here, ``$f \approx g$'' is used to indicate that the ratio of the two sides approaches one as $x$ tends to infinity, that is,
\[
\lim_{x\go\infty}\frac{\Gamma{(\gamma_1x +\gamma_2)}}{\sqrt{2\pi}e^{-\gamma_{1}x} (\gamma_{1}x )^{\gamma_{1}x + \gamma_{2} - 1/2}} =1.
\]

\begin{pf*}{Proof of Theorem \ref{schwarz:app}}
When the model dimension is $j$ is bounded and the sample size $n$ is large, it follows directly from Stirling's formula that
\[
\Gamma \biggl(\frac{n - j -1}{2} \biggr) \approx \sqrt{2\pi}e^{-n/2}
\biggl(\frac{n}{2} \biggr)^{(n - j)/2-1}\quad\mbox{and}\quad \Gamma
\biggl(\frac{n - 1}{2} \biggr) \approx \sqrt{2\pi}e^{-n/2} \biggl(
\frac{n}{2} \biggr)^{n/2-1}.
\]
The Bayes factor in (\ref{BFequation}) is asymptotically equivalent
\begin{eqnarray*}
\operatorname{BF}[M_j: M_i] &\approx& \frac{\sqrt{2\pi}e^{-n/2} (n/2 )^{(n - j)/2-1} }{\sqrt{2\pi}e^{-n/2} (n/2 )^{n/2-1}}
\bigl(1 - R_j^2 \bigr)^{-(n-j-1)/2+a+1}
\\
&\approx& \biggl(\frac{n}{2} \biggr)^{-j/2} \bigl(1 -
R_j^2 \bigr)^{-n/2}
\approx \exp \biggl[-\frac{j}{2} \log n - \frac{n}{2}\log \bigl(1 -
R_j^2 \bigr) \biggr].
\end{eqnarray*}
This completed the proof.
\end{pf*}

We now investigate the model selection consistency of Bayes factor in (\ref{nonnest:tw0}) under the three different asymptotic scenarios mentioned above. For simplicity of notation, let $c_i$ represent a finite constant for $i=1, 2, \ldots, 5$ throughout the following proofs. When $(j/2+a+1)$ and $(n-j-1)/2$ are sufficiently large, it follows directly from Stirling's formula that
\[
\Gamma \biggl(\frac{j}{2} + a + 1 \biggr) \approx \sqrt{2
\pi}e^{-j/2} \biggl(\frac{j}{2} \biggr)^{j/2 + a+1/2}
\]
and
\[
\Gamma \biggl(\frac{n - j -1}{2} \biggr) \approx \sqrt{2\pi}e^{-(n - j)/2}
\biggl(\frac{n - j}{2} \biggr)^{(n - j)/2-1}.
\]

\begin{pf*}{Proof of Theorem \ref{theorem:002}}
Under Scenario 1, $i = O(n^{a_1})$ and $j = O(n^{a_2})$ with $0 \leq a_1 \leq a_2 < 1$, by using the two approximation equations above, it follows that
\begin{eqnarray}
\label{thm:0201}
\nonumber
\operatorname{BF}[M_j: M_i] &=&
\frac{\Gamma (j/2+a+1 ) \Gamma ((n-j-1)/2 )}{\Gamma (i/2+a+1 ) \Gamma ((n-i-1)/2 )}\frac{ (1 - R_j^2 )^{-(n-j-1)/2+a+1}}{ (1 - R_i^2 )^{-(n-i-1)/2+a+1}}
\\
&=& c_1 \frac{j^{j/2+a+1}(n-j)^{(n-j)/2+1}}{i^{i/2+a+1}(n-i)^{(n-i)/2+1}}\frac{ (1 - R_j^2 )^{-(n-j)/2}}{ (1 - R_i^2 )^{-(n-i)/2}}
\\
&=& c_1 \frac{(j/n)^{j/2}}{(i/n)^{i/2}} \biggl(\frac{j}{i}
\biggr)^{a+1} \biggl(\frac{1-j/n}{1-i/n} \biggr) \biggl[\frac{(1-j/n)^{1-j/n}}{(1-i/n)^{1-i/n}}
\frac{ (1 - R_j^2 )^{-(1-j/n)}}{ (1 - R_i^2 )^{-(1-i/n)}} \biggr]^{n/2}.\nonumber
\end{eqnarray}
\begin{longlist}[(a)]
\item[(a)] We first show the Bayes factor consistency when the true model is $M_i$. As $n$ tends to infinity, we observe that the dominated term in brackets of equation (\ref{thm:0201}) can be approximated by
\begin{eqnarray*}
\frac{(1-j/n)^{1-j/n}}{(1-i/n)^{1-i/n}}\frac{ (1 - R_j^2 )^{-(1-j/n)}}{ (1 - R_i^2 )^{-(1-i/n)}} \approx \biggl(\frac{1 - R_j^2}{1 - R_i^2}
\biggr)^{-1},
\end{eqnarray*}
because of $j/n$ and $i/n$ approaching to zero as $n$ approaches infinity. From Lemma \ref{lemma:00002}(a) and the fact that $\delta_{ii}=0$, we observe that under $M_i$, it follows
\begin{eqnarray*}
\operatorname{BF}[M_j: M_i] &=& c_2
\frac{(j/n)^{j/2}}{(i/n)^{i/2}} \biggl(\frac{j}{i} \biggr)^{a+1} \biggl(
\frac{1-j/n}{1-i/n} \biggr) \biggl(\frac{1 + \delta_{ij}}{1 + \delta_{ii}} \biggr)^{-n/2}
\\
&=& c_2 \frac{(j/n)^{j/2}}{(i/n)^{i/2}} \biggl(\frac{j}{i}
\biggr)^{a+1} \biggl(\frac{1-j/n}{1-i/n} \biggr) (1 +
\delta_{ij} )^{-n/2},
\end{eqnarray*}
which approaches zero as $\delta_{ij} > 0$, indicating that the Bayes factor in (\ref{nonnest:tw0}) is consistent when $M_i$ is true.

\item[(b)] Consistency under $M_j$ is provided as follows. By using Lemma \ref{lemma:00002}(a), it follows that under model $M_j$, the Bayes factor in (\ref{nonnest:tw0}) can be further approximated by
\begin{eqnarray*}
\operatorname{BF}[M_j: M_i] &=& c_3
\frac{(j/n)^{j/2}}{(i/n)^{i/2}} \biggl(\frac{j}{i} \biggr)^{a+1} \biggl(
\frac{1-j/n}{1-i/n} \biggr) \biggl(\frac{1 + \delta_{jj}}{1 + \delta_{ji}} \biggr)^{-n/2}
\\
&=& c_3 \frac{(j/n)^{j/2}}{(i/n)^{i/2}} \biggl(\frac{j}{i}
\biggr)^{a+1} \biggl(\frac{1-j/n}{1-i/n} \biggr) (1 +
\delta_{ji} )^{n/2},
\end{eqnarray*}
because $\delta_{jj}=0$. It should be noted that as $n$ tends to infinity, the fifth dominated term approaches infinity if $\delta_{ji} > 0$. Therefore, the Bayes factor also approaches infinity as $\delta_{ji} >0$, proving the consistency under $M_j$. This completes the proof the theorem.\qed
\end{longlist}\noqed
\end{pf*}

\begin{pf*}{Proof of Theorem \ref{theorem:003}}
Under Scenario 2, $i = O(n^{a_1})$ and $j = O(n^{a_2})$ with $0 \leq a_1 < a_2 = 1$, by using the two approximation equations above, it follows that
\begin{eqnarray}
\label{thm:03001}
\nonumber
\operatorname{BF}[M_j: M_i] &=&
\frac{\Gamma (j/2+a+1 ) \Gamma ((n-j-1)/2 )}{\Gamma (i/2+a+1 ) \Gamma ((n-i-1)/2 )}\frac{ (1 - R_j^2 )^{-(n-j-1)/2+a+1}}{ (1 - R_i^2 )^{-(n-i-1)/2+a+1}}
\\
&=& c_4 \frac{(j/i)^{a+1}}{(i/n)^{i/2}} \biggl(\frac{1-j/n}{1-i/n} \biggr)
\\
&&{}\times  \biggl[
\biggl(\frac{j}{n} \biggr)^{j/n}\frac{(1-j/n)^{1-j/n}}{(1-i/n)^{1-i/n}}
\frac{ (1 - R_j^2 )^{-(1-j/n)}}{ (1 - R_i^2 )^{-(1-i/n)}} \biggr]^{n/2}.\nonumber
\end{eqnarray}
\begin{longlist}[(a)]
\item[(a)] If the true model is $M_i$, from Lemma \ref{lemma:00002}(b) and the fact that $\delta_{ii}=0$, we observe that the dominated term in brackets of (\ref{thm:03001}) can be approximated by
\begin{eqnarray*}
&& \biggl(\frac{j}{n} \biggr)^{j/n}\frac{(1-j/n)^{1-j/n}}{(1-i/n)^{1-i/n}}
\frac{ (1 - R_j^2 )^{-(1-j/n)}}{ (1 - R_i^2 )^{-(1-i/n)}}
\\
&&\quad \approx \biggl(\frac{1}{r} \biggr)^{1/r}
\biggl(1-\frac{1}{r} \biggr)^{1-1/r} \biggl(\frac{1 - R_j^2}{1 - R_i^2}
\biggr)^{-(1-1/r)} \bigl(1 - R_i^2
\bigr)^{1/r}
\\
&&\quad \approx \biggl(\frac{1}{r} \biggr)^{1/r} \biggl(
\frac{1-1/r}{1-1/r+\delta_{ij}} \biggr)^{1-1/r} \biggl(\frac{1}{1+\delta_{i0}}
\biggr)^{1/r}.
\end{eqnarray*}
Accordingly,  the approximation of Bayes factor in (\ref{nonnest:tw0}) is given by
\begin{eqnarray*}
\operatorname{BF}[M_j: M_i] &\approx& c_4
\frac{(j/i)^{a+1}}{(i/n)^{i/2}} \biggl[ \biggl(\frac{1}{r} \biggr)^{1/r}
\biggl(\frac{1-1/r}{1-1/r+\delta_{ij}} \biggr)^{1-1/r} \biggl(\frac{1}{1+\delta_{i0}}
\biggr)^{1/r} \biggr]^{n/2},
\end{eqnarray*}
which approaches zero as $n$ tends to infinity, and therefore, the consistency under $M_i$ is proved.

\item[(b)] If the true model is $M_j$, from Lemma \ref{lemma:00002}(b) and the fact that $\delta_{jj}=0$, we observe that the dominated term in brackets of (\ref{thm:03001}) can be approximated by
\begin{eqnarray*}
&& \biggl(\frac{j}{n} \biggr)^{j/n}\frac{(1-j/n)^{1-j/n}}{(1-i/n)^{1-i/n}}
\frac{ (1 - R_j^2 )^{-(1-j/n)}}{ (1 - R_i^2 )^{-(1-i/n)}}
\\
&&\quad \approx \biggl(\frac{1}{r} \biggr)^{1/r}
\biggl(1-\frac{1}{r} \biggr)^{1-1/r} \biggl(\frac{1 - R_j^2}{1 - R_i^2}
\biggr)^{-1} \bigl(1 - R_j^2
\bigr)^{1/r}
\\
&&\quad \approx \biggl(\frac{1}{r} \biggr)^{1/r} \biggl(1-
\frac{1}{r} \biggr)^{1-1/r} \biggl(\frac{1-1/r}{1+\delta_{ji}}
\biggr)^{-1} \biggl(\frac{1-1/r}{1+\delta_{j0}} \biggr)^{1/r}
\\
&&\quad \approx \biggl(\frac{1}{r} \biggr)^{1/r} (1 +
\delta_{ji}) \biggl(\frac{1}{1+\delta_{j0}} \biggr)^{1/r}.
\end{eqnarray*}
Therefore, the Bayes factor in (\ref{nonnest:tw0}) under $M_j$ turns out to be
\begin{eqnarray}
\label{thm:0203} \operatorname{BF}[M_j: M_i] &=&
c_5 \frac{(j/i)^{a+1}}{(i/n)^{i/2}} \biggl[ \biggl(\frac{1}{r}
\biggr)^{1/r} (1 + \delta_{ji}) \biggl(\frac{1}{1+\delta_{j0}}
\biggr)^{1/r} \biggr]^{n/2}.
\end{eqnarray}
To show the consistency under $M_j$, it is sufficient to show that the dominated term in brackets of~(\ref{thm:0203}) is strictly larger than one when $n$ tends to infinity. This is equivalent to
\[
\biggl(\frac{1}{r} \biggr)^{1/r} (1 + \delta_{ji})
\biggl(\frac{1}{1+\delta_{j0}} \biggr)^{1/r} > 1,
\]
which gives that
\[
\delta_{ji} > \bigl[r(1+\delta_{j0}) \bigr]^{1/r} -
1.
\]
On the other hand, we have $\delta_{ji} \leq \delta_{j0}$, which provides that
\[
\delta_{j0} \geq \delta_{ji} > \bigl[r(1+
\delta_{j0}) \bigr]^{1/r} - 1,
\]
indicating that
\[
\delta_{j0} > r^{1/(r-1)} - 1 = \delta(r).
\]
In order for the interval where the distance $\delta_{ji}$ should lie
\[
\delta_{ji} \in \bigl( \bigl[r(1+\delta_{j0})
\bigr]^{1/r} - 1, \delta_{j0} \bigr]
\]
to be nonempty, a necessary and sufficient condition is that $\delta_{j0} > \delta(r)$. This completes the proof.
\end{longlist}
\end{pf*}

\begin{pf*}{Proof of Theorem \ref{theorem:004}}
Under Scenario 3, $i = O(n^{a_1})$ and $j = O(n^{a_2})$ with $a_1 = a_2 = 1$, by using the two approximations equations, it follows that
\begin{eqnarray}
\label{thm:0301}
\nonumber
\operatorname{BF}[M_j: M_i] &=&
\frac{\Gamma (j/2+a+1 ) \Gamma ((n-j-1)/2 )}{\Gamma (i/2+a+1 ) \Gamma ((n-i-1)/2 )}\frac{ (1 - R_j^2 )^{-(n-j-1)/2+a+1}}{ (1 - R_i^2 )^{-(n-i-1)/2+a+1}}
\\
&=& c_5 \biggl(\frac{j}{i} \biggr)^{a+1} \biggl(
\frac{1-j/n}{1-i/n} \biggr)
\\
&&{}\times  \biggl[\frac{(j/n)^{j/n}}{(i/n)^{i/n}}\frac{(1-j/n)^{1-j/n}}{(1-i/n)^{1-i/n}}
\frac{ (1 - R_j^2 )^{-(1-j/n)}}{ (1 - R_i^2 )^{-(1-i/n)}} \biggr]^{n/2}.\nonumber
\end{eqnarray}
\begin{longlist}[(a)]
\item[(a)] If the true model is $M_i$, from Lemma \ref{lemma:00002}(c) and the fact that $\delta_{ii}=0$, we observe that the dominated term in brackets of (\ref{thm:0301}) can be approximated by
\begin{eqnarray}
\label{thm:0302}
\nonumber
&& \frac{(j/n)^{j/n}}{(i/n)^{i/n}}\frac{(1-j/n)^{1-j/n}}{(1-i/n)^{1-i/n}}\frac{ (1 - R_j^2 )^{-(1-j/n)}}{ (1 - R_i^2 )^{-(1-i/n)}}
\\
\nonumber
&&\quad \approx \frac{(1/r)^{1/r}}{(1/s)^{1/s}}\frac{(1-1/r)^{1-1/r}}{(1-1/s)^{1-1/s}} \biggl(
\frac{1 - R_j^2}{1 - R_i^2} \biggr)^{-(1-1/r)} \bigl(1 - R_i^2
\bigr)^{1/r-1/s}
\nonumber\\[-8pt]\\[-8pt]\nonumber
&&\quad \approx \frac{(1/r)^{1/r}}{(1/s)^{1/s}}\frac{(1-1/r)^{1-1/r}}{(1-1/s)^{1-1/s}} \biggl(
\frac{1 + \delta_{ij}-1/r}{1 - 1/s} \biggr)^{-(1-1/r)} \biggl(\frac{1 - 1/s}{1 + \delta_{i0}}
\biggr)^{1/r-1/s}
\\
&&\quad \approx \frac{(1/r)^{1/r}}{(1/s)^{1/s}}\frac{ [1 + \delta_{ij}/(1-1/r) ]^{-(1-1/r)}}{(1 +\delta_{i0})^{1/r-1/s}}.\nonumber
\end{eqnarray}
For the Bayes factor to be consistent, it is sufficient to show that the dominated term in (\ref{thm:0302}) is strictly less than 1 as $n$ approaches infinity. This is equivalent to
\[
\biggl(1 + \frac{\delta_{ij}}{1-1/r} \biggr)^{1-1/r} > \frac{(1/r)^{1/r}}{(1/s)^{1/s}} (1+
\delta_{i0} )^{1/s-1/r},
\]
which implies that
\[
\delta_{ij} > \frac{r-1}{r} \biggl\{ \biggl[\frac{s^{1/s}}{r^{1/r}}(1+
\delta_{i0})^{1/s - 1/r} \biggr]^{r/(r-1)} - 1 \biggr\}.
\]
In addition, from the property of the pseudo-distance, we have $\delta_{i0} \geq \delta_{ij}$. Therefore, it follows that
\[
\delta_{i0} \geq \delta_{ij} > \frac{r-1}{r} \biggl\{
\biggl[\frac{s^{1/s}}{r^{1/r}}(1+\delta_{i0})^{1/s - 1/r}
\biggr]^{r/(r-1)} - 1 \biggr\},
\]
indicating that the value of $\delta_{ij}$ must satisfy
\[
\biggl(1 + \frac{\delta_{i0}}{1-1/r} \biggr)^{1-1/r} > \frac{(1/r)^{1/r}}{(1/s)^{1/s}} (1+
\delta_{i0} )^{1/s-1/r}.
\]
Under the conditions stated in the theorem, we take limits and obtain that the Bayes factor tends to zero, and thus, the Bayes factor is consistent under $M_i$.
\item[(b)] If the true model is $M_j$, from Lemma \ref{lemma:00002}(c) and the fact that $\delta_{jj}=0$, we observe that the dominated term in brackets of (\ref{thm:0301}) can be approximated by
\begin{eqnarray}
\label{thm:0303}
\nonumber
&& \frac{(j/n)^{j/n}}{(i/n)^{i/n}}\frac{(1-j/n)^{1-j/n}}{(1-i/n)^{1-i/n}}
\frac{ (1 - R_j^2 )^{-(1-j/n)}}{ (1 - R_i^2 )^{-(1-i/n)}}
\\
\nonumber
&&\quad \approx \frac{(1/r)^{1/r}}{(1/s)^{1/s}}\frac{(1-1/r)^{1-1/r}}{(1-1/s)^{1-1/s}} \biggl(
\frac{1 - R_j^2}{1 - R_i^2} \biggr)^{-(1-1/s)} \bigl(1 - R_j^2
\bigr)^{1/r-1/s}
\nonumber\\[-8pt]\\[-8pt]\nonumber
&&\quad \approx \frac{(1/r)^{1/r}}{(1/s)^{1/s}}\frac{(1-1/r)^{1-1/r}}{(1-1/s)^{1-1/s}} \biggl(
\frac{1 -1/r}{1 + \delta_{ji} - 1/s} \biggr)^{-(1-1/s)} \biggl(\frac{1 - 1/r}{1 + \delta_{j0}}
\biggr)^{1/r-1/s}
\\
&&\quad \approx \frac{(1/r)^{1/r}}{(1/s)^{1/s}}\frac{ [1+\delta_{ji}/(1-1/s) ]^{1-1/s}}{(1+\delta_{j0})^{1/r-1/s}}.\nonumber
\end{eqnarray}
For the Bayes factor to be consistent, it is sufficient to show that the dominated term in (\ref{thm:0303}) is strictly larger than one as $n$ approaches infinity. This is equivalent to
\[
\frac{(1/r)^{1/r}}{(1/s)^{1/s}}\frac{ [1+\delta_{ji}/(1-1/s) ]^{1-1/s}}{(1+\delta_{j0})^{1/r-1/s}} > 1.
\]
Simple algebra shows that
\[
\delta_{ji} > \frac{s-1}{s} \biggl\{ \biggl[\frac{r^{1/r}}{s^{1/s}} (1
+ \delta_{j0} )^{1/r-1/s} \biggr]^{s/(s-1)}-1 \biggr\}.
\]
On the other hand, we also have $\delta_{j0} \geq \delta_{ji}$, which provides that
\begin{equation}
\label{th04:02} \delta_{j0} \geq \delta_{ji} >
\frac{s-1}{s} \biggl\{ \biggl[\frac{r^{1/r}}{s^{1/s}} (1 + \delta_{j0}
)^{1/r-1/s} \biggr]^{s/(s-1)}-1 \biggr\},
\end{equation}
indicating that
\[
\biggl(1 + \frac{\delta_{j0}}{1-1/s} \biggr)^{1-1/s}> \frac{r^{1/r}}{s^{1/s}}(1+
\delta_{j0})^{1/r-1/s}.
\]
In order for the interval where the distance $\delta_{ji}$ should lie
\[
\delta_{ji} \in \biggl(\frac{s-1}{s} \biggl[\frac{r^{1/r}}{s^{1/s}} (1 +
\delta_{j0} )^{1/r-1/s} -1 \biggr]^{s/(s-1)},
\delta_{j0} \biggr]
\]
to be nonempty, a necessary and sufficient condition is that $\delta_{j0}$ satisfies inequality in (\ref{th04:02}). This completes the proof.\qed
\end{longlist}\noqed
\end{pf*}
\end{appendix}

\section*{Acknowledgements}
The authors thank the Editor and two referees for their helpful comments, which have led to an improvement of the manuscript. 



\printhistory
\end{document}